\documentclass[11pt]{amsart}
\usepackage[margin=3cm]{geometry} 
\usepackage{amsmath}
\usepackage{subcaption}
\usepackage{amssymb}
\usepackage{amsthm}
\usepackage{float}
\usepackage{ulem}
\usepackage{xcolor}
\usepackage[utf8x]{inputenc}
\usepackage{comment}
\usepackage{framed}
\usepackage{caption}

\newcommand{\eps}{\varepsilon}

\DeclareMathAlphabet{\mathbbb}{U}{bbold}{m}{n}

\newtheorem{thm}{Theorem}[section]
\newtheorem{cor}[thm]{Corollary}

\newtheorem{lem}[thm]{Lemma}

\theoremstyle{definition}

\theoremstyle{remark}
\newtheorem{rem}[thm]{Remark}

\numberwithin{equation}{section}
\usepackage{graphicx}

\usepackage{enumerate}

\newcommand{\matrice}{\begin{pmatrix}}
\newcommand{\ok}{\end{pmatrix}}

\begin{document}
\title[]{Remarks on Weyl-type bounds for Steklov eigenvalues}

\dedicatory{Dedicated to Alessandro Savo on the occasion of his 70th birthday}

\author[Provenzano]{Luigi Provenzano}
\address{Dipartimento di Scienze di Base e Applicate per l'Ingegneria, Sapienza Universit\`a di Roma, Via Scarpa 12 - 00161 Roma, Italy, e-mail: {\sf luigi.provenzano@uniroma1.it}.}

\begin{abstract}
We prove that, for $n\geq 3$, there is no constant $C_n>0$ depending only on $n$ such that the Steklov eigenvalues $\sigma_k(\Omega)$ satisfy
$|\partial\Omega|^{\frac 1{n-1}}\sigma_k(\Omega)\leq C_nk^{\frac1{n-1}}$ for every smooth bounded domain $\Omega\subset\mathbb R^n$ and every $k\geq1 $, providing a negative answer to an open problem posed by Girouard and Polterovich \cite{GiPo}. On the other hand, we prove that there exists a constant $C_n>0$ depending only on $n$ such that $|\Omega|^{\frac 1n}\sigma_k(\Omega)\leq C_nk^{\frac1{n-1}}$ for every smooth bounded domain $\Omega\subset\mathbb R^n$ and every $k\geq 1$. This estimate, combined with the isoperimetric bound of Colbois, El Soufi and Girouard \cite{colboisgirouard_steklov}, implies the bound 
$|\partial\Omega|^{\frac 1{n-1}}\sigma_k(\Omega)\leq C_nk^{\frac1{n-1}+\frac{n-2}{n(n-1)^2}}$, where the exponent of $k$ turns out to be sharp.
\end{abstract}


\thanks{The author acknowledges the support of the INdAM GNSAGA group.}

\maketitle

\section{Introduction and description of the  results}

In this note we discuss Weyl-type upper bounds for Steklov eigenvalues. Throughout the note, $\Omega$ will be a smooth, bounded domain in $\mathbb R^n$, $n\geq 2$. The Steklov eigenvalue problem on $\Omega$ is:
\begin{equation}\label{steklov}
\begin{cases}
\Delta u=0 \,, & {\rm in\ }\Omega\,,\\
\partial_Nu=\sigma u\,, & {\rm on\ }\partial\Omega.
\end{cases}
\end{equation}
Here $\Delta=-\sum_{i=1}^n\partial_{ii}^2$ and $N$ is the unit outer normal to $\Omega$. Problem \eqref{steklov} admits an increasing sequence of non-negative eigenvalues of finite multiplicity:
$$
0=\sigma_1(\Omega)<\sigma_2(\Omega)\leq\cdots\leq\sigma_k(\Omega)\leq\cdots\nearrow+\infty.
$$
The eigenfunctions can be chosen so that their traces form an orthonormal basis of $L^2(\partial\Omega)$. We refer to \cite{CGGS,GiPo} for a panoramic view on the Steklov problem.

\smallskip

It is well-known that Steklov eigenvalues satisfy the following Weyl's law (see e.g., \cite{GiPo}):
\begin{equation}\label{weyl}
\lim_{k\to\infty}\sigma_k(\Omega)\left(\frac{|\partial\Omega|}{k}\right)^{\frac{1}{n-1}}=2\pi\omega_{n-1}^{-\frac{1}{n-1}}, 
\end{equation}
which can be re-written:
$$
|\partial\Omega|^\frac{1}{n-1}\sigma_k(\Omega)=2\pi\omega_{n-1}^{-\frac{1}{n-1}}k^{\frac{1}{n-1}}+o(k^{\frac{1}{n-1}})\,,\ \ \ k\to\infty.
$$
Here $\omega_n$ denotes the volume of the unit ball in $\mathbb R^n$. Steklov eigenvalues scale as follows under homotheties: $\sigma_k(\gamma\Omega)=\gamma^{-1}\sigma_k(\Omega)$ for $\gamma>0$. Therefore, the natural scale-invariant quantity is $|\partial\Omega|^{\frac{1}{n-1}}\sigma_k(\Omega)$, which is usually referred to as (perimeter-){\it normalised eigenvalue}. The natural question is then if one can bound from above $|\partial\Omega|^{\frac{1}{n-1}}\sigma_k(\Omega)$ in terms of $n$ and $k$ alone.

\smallskip

The literature on upper bounds for Steklov eigenvalues is quite vast, see, among others, \cite{CoElGi,CoGiAIF,ColGirGit,CGH,ColGit,GiPosurfaces,hassannezhadJFA,Prov_Stub}. We refer to the survey articles \cite{CGGS} and \cite{GiPo} for an exhaustive treatment. We recall just a few pivotal results, confining the (brief) discussion to the case of Euclidean domains. Colbois, El Soufi and Girouard \cite{colboisgirouard_steklov} proved the following upper bound (see also Hassannezhad \cite{hassannezhadJFA} for related conformal upper bounds):
\begin{equation}\label{isoperimetric}
|\partial\Omega|^{\frac{1}{n-1}}\sigma_k(\Omega)\leq C_n k^{\frac 2n}.
\end{equation}
Hence the normalised eigenvalues are uniformly bounded from above in terms of $n$ and $k$ only. In dimension $2$, inequality \eqref{isoperimetric} states that there exists a universal constant $C$ such that $|\partial\Omega|\sigma_k(\Omega)\leq Ck$, and this bound is compatible with the Weyl's law \eqref{weyl}. For $n\geq 3$, the exponent $\frac 2n$ in \eqref{isoperimetric} is worse than the exponent $\frac {1}{n-1}$ in \eqref{weyl}. A natural question is:

\smallskip
{\it does there exist a constant $C_n$ depending only on $n$ such that any smooth, bounded Euclidean domain $\Omega$ satisfies}
\begin{equation}\label{conj}
|\partial\Omega|^{\frac{1}{n-1}}\sigma_k(\Omega)\leq C_nk^{\frac{1}{n-1}}{\rm\ ?}
\end{equation}
\smallskip

This is exactly Open Problem 5 in \cite{GiPo} and also Remark 1.11 in \cite{GiLa}. Note that \eqref{conj} is not a question about the asymptotic behavior of the Steklov eigenvalues of a single domain, but rather about the existence of an upper estimate compatible with Weyl's law \eqref{weyl} that is uniform over all Euclidean domains and indices $k$. Thus, any counterexample should involve a family of domains degenerating (in some sense) with $k$. The question is specific to dimension $n\geq 3$ since, as already observed, inequality \eqref{isoperimetric} proved in \cite{colboisgirouard_steklov} coincides with \eqref{conj} for $n=2$.

\smallskip

A first partial answer was given in a paper by Stubbe and the author \cite{Prov_Stub}. The inequality is stated as follows:
\begin{equation}\label{PS}
\sigma_k(\Omega)\leq C_n\left(\frac{k}{|\partial\Omega|}\right)^{\frac{1}{n-1}}+K_\Omega,
\end{equation}
where $K_\Omega>0$ depends in an explicit way on the geometry of $\Omega$ near a tubular neighborhood of $\partial\Omega$ through the mean curvature and the rolling radius. Precisely, $K_\Omega={\rm roll}(\Omega)^{-1}+n\|\mathcal H\|_{L^\infty(\partial\Omega)}$, where ${\rm roll}(\Omega)$ is the rolling radius of $\Omega$ and $\mathcal H$ is the mean curvature of $\partial\Omega$. Inequality \eqref{PS} follows from a comparison of the Steklov eigenvalues of $\Omega$ and the (square roots of the) Laplace eigenvalues of $\partial\Omega$. Hence a Weyl-type upper bound holds, and in particular, the Weyl-type behavior in $k$ of the upper bound \eqref{PS} is separated from the geometry of the domain through an explicit {\it additive constant}. This bound has been extended in the context of Riemannian manifolds in \cite{CGH,xiong}. In this more general setting, the additive constant depends on the principal curvatures of the boundary, the rolling radius and on bounds on the sectional curvatures of the manifold in a neighborhood of the boundary. It is worth mentioning that in the context of Riemannian manifolds the validity of a bound of the form \eqref{conj} cannot hold: a geometric additive constant as in \eqref{PS} is indeed needed if $n\geq 3$. This is proved by Colbois, El Soufi and Girouard in  \cite{CoElGi}, where they prove the existence of conformal metrics on compact manifolds with boundary (which coincide with a fixed reference metric on the boundary) with arbitrarily large first non-trivial Steklov eigenvalue. However, this does not settle the problem in $\mathbb R^n$, $n\geq 3$, where the question remained open.

\smallskip

 In this note we prove that no inequality of the form \eqref{conj} can hold for Euclidean domains. Here is the first main theorem:

\begin{thm}\label{main}
Let $n\geq 3$. Then there exists a sequence of smooth, bounded domains $\{\Omega_k\}_{k=1}^{\infty}$ in $\mathbb R^n$ such that, for every $k\geq 2$
\begin{equation}\label{lower_b}
|\partial\Omega_k|^{\frac{1}{n-1}}\sigma_k(\Omega_k)\geq C_n k^{\frac{1}{n-1}+\frac{n-2}{n(n-1)^2}},
\end{equation}
where $C_n>0$ depends only on $n$. In particular,
$$
\lim_{k\to\infty}\sigma_k(\Omega_k)\left(\frac{|\partial\Omega_k|}{k}\right)^{\frac{1}{n-1}}=+\infty.
$$
\end{thm}

Theorem \ref{main} is a consequence of a homogenisation process linking Steklov and Neumann eigenvalues, which has been discovered by Girouard, Henrot and Lagacé in \cite{GHL}. The domains $\Omega_k$ are obtained by periodically perforating a fixed domain (the unit ball $B$) with identical spherical holes whose centers are placed along a grid of size $\eps_k$ (sufficiently small) and with radii $r_k=b_n k^{\frac{n-2}{n(n-1)^2}}
\eps_k^{\frac{n}{n-1}}$, where $b_n>0$ is some explicit dimensional constant.

\smallskip

Note that Theorem \ref{main} does not violate \eqref{PS}. On the contrary, it shows that the additive geometric constant cannot be removed. In fact, if $\mathcal H_k$ denotes the mean curvature of  $\partial\Omega_k$ and ${\rm roll}(\Omega_k)$ the rolling radius, we have that $\|\mathcal H_k\|_{L^\infty(\partial\Omega_k)}\to\infty$ and ${\rm roll}(\Omega_k)\to 0$ as $k\to\infty$.

\smallskip It is not surprising that it is the homogenisation process introduced in \cite{GHL} that allows one to produce the family of domains $\{\Omega_k\}_{k=1}^\infty$ failing \eqref{conj}. In \cite{GHL} the homogenisation process is used to link Steklov and Neumann eigenvalues through an intermediate problem, called {\it dynamical problem} (see Problem \eqref{dynamical}). This link is used to transfer isoperimetric inequalities and upper bounds from Steklov to Neumann eigenvalues, and in particular to produce Euclidean domains with {\it large} normalised Steklov eigenvalues. This technique was then extended to manifolds in \cite{GiLa}, where many striking consequences on large Steklov eigenvalues, on free-boundary minimal hypersurfaces, and on the sharpness of known upper bounds for Steklov eigenvalues are deduced. This technique has later been included in a more general framework of continuity of eigenvalues depending on Radon measures \cite{GiKaLa}.

\smallskip

The argument used in the proof of Theorem \ref{main} exploits the homogenisation theorem of Girouard, Henrot and Lagacé in a different way. In \cite{GHL}, as the homogenisation parameter tends to zero, the Steklov eigenvalues of a periodically perforated domain converge to those of the dynamical problem \eqref{dynamical}, where the parameter $\beta>0$ is determined by the relative size of the holes. Moreover, for each fixed $k$, as $\beta\to\infty$,  the $k$-th dynamical eigenvalue rescaled by $\beta$ converges to the corresponding Neumann eigenvalue. In our construction, rather than taking this limit with $k$ fixed, we exploit, in the dynamical problem \eqref{dynamical}, the competition between the boundary and interior contributions by allowing $\beta$ to depend on $k$. More precisely, we choose $\beta_k\sim k^{\frac{n-2}{n(n-1)}}$. Here, for a positive sequence $s_k$, we write $s_k\sim k^a$ if there exist positive constants $c_{n,1},c_{n,2}$, depending only on $n$, such that $c_{n,1}k^a\leq s_k\leq c_{n,2}k^a$. With this choice, the $k$-th dynamical eigenvalue is bounded from below by $a_n k^{\frac1{n-1}}$, where $a_n>0$ is a dimensional constant. For each sufficiently large $k$, we then choose the homogenisation parameter $\eps_k$ sufficiently small. The resulting domain $\Omega_k$ satisfies $|\partial\Omega_k|\sim\beta_k\sim k^{\frac{n-2}{n(n-1)}}$. Putting everything together, we obtain the lower bound \eqref{lower_b}. The estimates for the dynamical eigenvalues are obtained by linking the dynamical spectrum to a suitable Robin spectrum.

\smallskip

We now discuss the sharpness of \eqref{lower_b}. The estimate \eqref{isoperimetric} of Colbois, El Soufi and Girouard is a consequence of the following isoperimetric bound:
\begin{equation}\label{isoperimetric2}
|\partial\Omega|^{\frac{1}{n-1}}\sigma_k(\Omega)\leq\frac{C_n}{I(\Omega)^{\frac{n-2}{n-1}}}k^{\frac 2n},
\end{equation}
where $I(\Omega)$ is the {\it isoperimetric ratio} of $\Omega$: $I(\Omega)=\frac{|\partial\Omega|}{|\Omega|^{\frac{n-1}{n}}}$. In our example $|\Omega_k|\sim 1$ and $|\partial\Omega_k|\sim k^{\frac{n-2}{n(n-1)}}$, and consequently the right hand side of \eqref{isoperimetric2} is $\sim k^{\frac{1}{n-1}+\frac{n-2}{n(n-1)^2}}$. Combining with \eqref{lower_b} we see that, along the sequence $\Omega_k$,
$$
|\partial\Omega_k|^{\frac 1{n-1}}\sigma_k(\Omega_k)\sim k^{\frac{1}{n-1}+\frac{n-2}{n(n-1)^2}}.
$$
Thus, the family of domains $\Omega_k$ in Theorem \ref{main} saturates the combined dependence on $k$ and $I(\Omega)$ of the isoperimetric bound \eqref{isoperimetric2}. The exponent $\frac 2n$ in the full isoperimetric estimate \eqref{isoperimetric2}, as well as the exponents of the volume and boundary measure, was proved to be sharp by Girouard and Lagacé \cite[Corollary 1.10]{GiLa} in the context of compact Riemannian manifolds. The construction in Theorem \ref{main} shows that the exponents in \eqref{isoperimetric2} are sharp also in the context of Euclidean domains.

\smallskip

This observation, however, does not rule out  the possibility that the exponent $\frac{2}{n}$ in \eqref{isoperimetric} is optimal for Euclidean domains, and this leads to the second question considered in this note:

\smallskip

{\it is the exponent $\frac{2}{n}$ in the right-hand side of \eqref{isoperimetric} sharp for Euclidean domains? If not, is $\frac{1}{n-1}+\frac{n-2}{n(n-1)^2}$ the sharp exponent?}

\smallskip

We prove that the sharp exponent is in fact $\frac{1}{n-1}+\frac{n-2}{n(n-1)^2}$. This is a consequence of our second main result:

\begin{thm}\label{main2}
There exists a constant $C_n>0$ depending only on $n$ such that, for every smooth, bounded domain of $\mathbb R^n$ and every $k\geq 1$
\begin{equation}\label{bound_true}
|\Omega|^\frac 1n\sigma_k(\Omega)\leq C_n k^{\frac 1{n-1}}.
\end{equation}
\end{thm}
Theorem \ref{main2} is an upper bound for the {\it volume}-normalised Steklov eigenvalues: the quantity $|\Omega|^{\frac{1}{n}}\sigma_k(\Omega)$ is in fact scaling invariant. The bound \eqref{bound_true} is weaker than the conjectured (and, as we have seen, false) bound \eqref{conj}: the isoperimetric inequality tells us that $|\partial\Omega|^{1/(n-1)}\geq n^{1/(n-1)}\omega_n^{1/n(n-1)}|\Omega|^{1/n}$. An immediate corollary of Theorem \ref{main2} is the following:

\begin{cor}\label{cor_isop_ratio}
There exists a constant $C_n>0$ depending only on $n$ such that, for every smooth, bounded domain of $\mathbb R^n$ and every $k\geq 1$
\begin{equation}\label{bound_true_2}
|\partial\Omega|^{\frac{1}{n-1}}\sigma_k(\Omega)\leq C_n I(\Omega)^{\frac{1}{n-1}}k^{\frac 1{n-1}}.
\end{equation}
\end{cor}
This inequality has to be compared with \eqref{isoperimetric2}, which is rewritten as follows:
\begin{equation}\label{comparison1}
|\partial\Omega|^{\frac{1}{n-1}}\sigma_k(\Omega)\leq C_n I(\Omega)^{-\frac{n-2}{n-1}}k^{\frac 2n}.
\end{equation}
Here we are denoting by $C_n$ a positive constant which depends only on $n$ and which may have been re-defined. In the isoperimetric bound \eqref{comparison1}, a larger isoperimetric ratio improves the bound, while in  \eqref{bound_true_2}, a larger isoperimetric ratio worsens the bound. The right-hand sides of \eqref{bound_true_2}-\eqref{comparison1} coincide (up to dimensional constants) when $I(\Omega)\sim k^\frac{n-2}{n(n-1)}$. This implies the following:

\begin{cor}\label{corollary_exponent}
There exists a constant $C_n>0$ such that, for every smooth, bounded domain in $\mathbb R^n$ and every $k\geq 1$
\begin{equation}\label{true_exponent}
|\partial\Omega|^{\frac{1}{n-1}}\sigma_k(\Omega)\leq C_n k^{\frac{1}{n-1}+\frac{n-2}{n(n-1)^2}}.
\end{equation}
\end{cor}
Note that Theorem \ref{main} shows that the exponent $\frac{1}{n-1}+\frac{n-2}{n(n-1)^2}$ in \eqref{true_exponent} is sharp. Moreover, inspecting the proof of Theorem \ref{main}, we observe that $\sigma_k(\Omega_k)\geq c_n k^\frac{1}{n-1}$, while $|\Omega_k|\sim 1$, hence also the exponent $\frac{1}{n-1}$ in \eqref{bound_true} and \eqref{bound_true_2} is sharp.

\smallskip

The upper bound \eqref{true_exponent} should be compared with \eqref{isoperimetric}, whose exponent $\frac2n$ was
proved to be sharp for compact Riemannian manifolds in \cite{GiLa}. Let us
stress that in the Riemannian setting no bound of the form \eqref{bound_true} can hold: applying
the homogenisation theorem of \cite{GiLa} e.g., to the round sphere $\mathbb S^n$ with constant weight
$\beta\equiv1$ produces domains $\Omega_{\varepsilon}$ with
$|\Omega_{\varepsilon}|,|\partial\Omega_{\varepsilon}|\to|\mathbb S^n|$ and
$\sigma_k(\Omega_{\varepsilon})\to\lambda_k(\mathbb S^n)$, where $\lambda_k(\mathbb S^n)$ are the Laplace eigenvalues on $\mathbb S^n$; then Weyl's law gives
$|\Omega_{\eps_k}|^{1/n}\sigma_k(\Omega_{\eps_k})\sim k^{\frac2n}$ along a suitable diagonal sequence. The reason is
that the boundary of $\Omega_{\varepsilon}$ consists of the boundaries of the removed
balls only: at a macroscopic scale the perforations disappear and the whole boundary
measure is redistributed in the bulk of $\mathbb S^n$. This cannot happen for a bounded Euclidean domain $\Omega$: indeed, after normalising $|\Omega|=1$, in the regime $k^a\geq |\partial\Omega|$, $0<a\leq \frac{1}{n-1}$, the presence of an ``exterior'' (the unbounded connected component of $\overline{\Omega^c}$) forces the existence of at least $k$ well-separated cubes of side $r\sim k^{-\frac{1}{n-1}}$ in which both the domain and its complement occupy a definite fraction of the volume of the cube. By the relative isoperimetric inequality, each such cube contains at least $c_n/k$ of boundary measure. This geometric fact, quantified in Lemma \ref{lem_cap}, is the mechanism behind Theorem \ref{main2}. A closed manifold carries no such ``outer boundary''; see Remark \ref{rem:riemannian} for a discussion of where the argument breaks down.


\smallskip

By scaling, it is sufficient to prove Theorem \ref{main2}  for domains $\Omega$ with $|\Omega|=1$. If $|\partial\Omega|\geq k^{\frac{n-2}{n(n-1)}}$, the statement of the theorem is a consequence of \eqref{isoperimetric2}. If $|\partial\Omega|\leq k^{\frac{n-2}{n(n-1)}}$, Lemma \ref{lem_cap} allows to find $k$ cubes $Q_i$ of side $\sim k^{-\frac{1}{n-1}}$, pairwise disjoint and such that $|\partial\Omega\cap Q_i|\geq c_n/k$ for some positive dimensional constant $c_n$. Then to each cube we associate a test function for the Rayleigh quotient of $\sigma_k(\Omega)$ which can be estimated explicitly in terms of $n$ and $k$ only.


\smallskip 

The present note is organised as follows: in Section \ref{proof} we present the proof of Theorem \ref{main} and in Section \ref{proof_main2} we present the proof of Theorem \ref{main2}. In Section \ref{23} we prove Lemmas \ref{lem_comparison}, \ref{countingR}, \ref{lowerR} (used in Section \ref{proof}),  and Lemma \ref{lem_cap} (used in Section \ref{proof_main2}).

\section{Proof of Theorem \ref{main}}\label{proof}

Throughout this section we assume $n\geq 3$.
\subsection{Construction of the domains}
The domains in Theorem \ref{main} are obtained by periodically perforating a fixed domain. Without loss of generality we take the fixed domain to be the unit ball $B:=B(0,1)$. Here $B(x_0,r)$ denotes the ball of radius $r$ in $\mathbb R^n$ centered in $x_0$. The construction is explained in \cite{GHL}, and we briefly repeat it here for the reader's convenience.

\smallskip Let $\eps\in(0,1)$, $m\in\mathbb Z^n$, and consider the cubes:
$$
Q_m^\eps:=\eps m+[-\eps/2,\eps/2]^n.
$$
For each $\eps\in(0,1)$, consider the set $I^\eps:=\{m\in\mathbb Z^n:Q_m^\eps\subset B\}$. Let $\beta>0$ and set
$$
r_{\beta,\eps}:=\beta^{1/(n-1)}\eps^{n/(n-1)}.
$$
With this choice
$$
r_{\beta,\eps}^{n-1}\eps^{-n}=\beta
$$
for all $\eps\in(0,1)$. Now, for any $\eps\in(0,\eps_0(\beta))$ we have that $r_{\beta,\eps}<\eps/2$, where $\eps_0(\beta):=\min\{1,\beta^{-1}2^{1-n}\}$. Define the `spherical holes' centered at $\eps m$ as
$$
\omega_m^{\beta,\eps}:=B(\eps m,r_{\beta,\eps}).
$$
Finally, define for all $\eps\in(0,\eps_0(\beta))$:
\begin{equation}\label{domains_eps}
\Omega_{\beta,\eps}:=B\setminus\overline{\bigcup_{m\in I^\eps}\omega_m^{\beta,\eps}}.
\end{equation}
The domains $\Omega_k$ of Theorem \ref{main} are domains of the form $\Omega_{\beta_k,\eps_k}$, with suitable choices of $\beta_k,\eps_k$.

\subsection{homogenisation process: the result of Girouard-Henrot-Lagacé}

We state here the main result of \cite{GHL}, which characterizes the asymptotic behavior of the Steklov problem \eqref{steklov} on $\Omega_{\beta,\eps}$ as $\eps\to 0$. To do so, we first introduce the {\it dynamical eigenvalue problem} (see \cite{F-VB}):
\begin{equation}\label{dynamical}
\begin{cases}
\Delta U=S_n\beta\Sigma U\,, & {\rm in\ }\Omega\,,\\
\partial_NU=\Sigma U\,, & {\rm on\ }\partial\Omega.
\end{cases}
\end{equation}
Here $\beta>0$ and $S_n=n\omega_n$ is the $n-1$-dimensional volume of the unit sphere $\mathbb S^{n-1}=\partial B$ in $\mathbb R^n$. Note that here the eigenvalue $\Sigma$ appears also in the boundary condition. When $\Omega$ is a smooth, bounded domain, Problem \eqref{dynamical} admits an increasing sequence of non-negative eigenvalues of  finite multiplicity:
$$
0= \Sigma_{1,\beta}(\Omega)<\Sigma_{2,\beta}(\Omega)\leq\cdots\leq\Sigma_{k,\beta}(\Omega)\leq\cdots\nearrow+\infty.
$$
We are now ready to state the main result in \cite{GHL}  (Theorem 2), applied to our specific situation.
\begin{thm}
Let $\beta>0$ and let $\Omega_{\beta,\eps}$ be defined by \eqref{domains_eps}. Then, for all $k\in\mathbb N$, $k\geq 1$:
$$
\lim_{\eps\to 0}\sigma_k(\Omega_{\beta,\eps})=\Sigma_{k,\beta}(B).
$$
\end{thm}

\begin{rem}In \cite{GHL} the fixed domain is not necessarily a ball; moreover, the authors also discuss convergence of eigenfunctions, which we don't use in this note.
\end{rem}

\subsection{The dynamical and the Robin problems}
 In this subsection we relate the dynamical eigenvalues to those of a Robin problem, and state some useful lemmas on lower bounds for Robin eigenvalues of the ball.

\smallskip

We recall the Robin eigenvalue problem on a smooth, bounded domain $\Omega$ of $\mathbb R^n$:
\begin{equation}\label{robin}
\begin{cases}
\Delta u=\lambda u \,, & {\rm in\ }\Omega\,,\\
\partial_Nu=\alpha u\,, & {\rm on\ }\partial\Omega,
\end{cases}
\end{equation}
where $\alpha\in\mathbb R$. In this note we will consider the case $\alpha\geq 0$.

The spectrum is discrete, and consists of an increasing sequence of eigenvalues of finite multiplicity, bounded from below:
$$
-\infty<\lambda_{1,\alpha}(\Omega)<\lambda_{2,\alpha}(\Omega)\leq\cdots\leq\lambda_{k,\alpha}(\Omega)\leq\cdots\nearrow+\infty.
$$

We state a few useful lemmas that we will employ to conclude the proof of Theorem \ref{main}.

\begin{lem}\label{lem_comparison}
Let $\Omega$ be a smooth, bounded domain in $\mathbb R^n$. Let $\beta>0$. Then 
$$
S_n\beta\Sigma_{k,\beta}(\Omega)=\lambda_{k,\Sigma_{k,\beta}(\Omega)}(\Omega).
$$
\end{lem}

\begin{lem}\label{countingR}
Let $N^R_\alpha(E):=\#\{\lambda_{j,\alpha}(B)<E\}$ be the Robin counting function. Then, for all $\alpha,E\geq 0$ we have
\begin{equation}\label{counting_R}
N^R_\alpha(E)\leq C_n(1+E^{n/2}+(\alpha^2+E)^{(n-1)/2})
\end{equation}
where $C_n>0$ depends only on $n$.
\end{lem}
A consequence of Lemma \ref{countingR} is the following:
\begin{lem}\label{lowerR}
Let $C_n$ be the constant appearing in \eqref{counting_R}. Let
\begin{eqnarray*}
&&\rho_n:=(8(1+C_n))^{-2/n}\,,\\
&&\eta_n:=\frac{1}{2}(8(1+C_n))^{-1/(n-1)}\,,\\
&&k_0(n):=\left[2^{4n(n-1)}(1+C_n)\right]+1\,,
\end{eqnarray*}
where $[\cdot]$ denotes the integer part of a real number.
Define, for any $k\in\mathbb N$, $k\geq 1$:
\begin{equation*}
\alpha_k:=\eta_n k^{1/(n-1)}.
\end{equation*}
Then, for any $k\geq k_0(n)$:
\begin{equation}\label{lower_R}
\lambda_{k,\alpha_k}(B)\geq \rho_n k^{2/n}.
\end{equation}
\end{lem}

\subsection{Conclusion of the proof of Theorem \ref{main}}

For any $k\in\mathbb N$, $k\geq k_0(n)$, define
$$
\beta_k:=\frac{\rho_n}{S_n\eta_n}k^{\frac{n-2}{n (n-1)}}=\frac{\rho_nk^{2/n}}{S_n\alpha_k},
$$
where $k_0(n),\rho_n,\eta_n,\alpha_k$ are defined in Lemma \ref{lowerR}. From Lemma \ref{lem_comparison} we have:
$$
\Sigma_{k,\beta_k}(B)=\frac{1}{S_n\beta_k}\lambda_{k,\Sigma_{k,\beta_k}(B)}(B).
$$
We claim that
\begin{equation}\label{final_claim}
\Sigma_{k,\beta_k}(B)\geq\alpha_k=\eta_nk^{1/(n-1)}.
\end{equation}
The validity of \eqref{final_claim} allows to conclude the proof of Theorem \ref{main}. In fact, from \cite[Formula (88)]{GHL} we have
$$
\lim_{\eps\to 0}|\partial\Omega_{\beta_k,\eps}|=|\partial B|+\frac{\rho_n}{\eta_n}k^{\frac{n-2}{n(n-1)}}|B|,
$$
which implies
$$
\lim_{\eps\to 0}\sigma_k(\Omega_{\beta_k,\eps})\left(\frac{|\partial\Omega_ {\beta_k,\eps}|}{k}\right)^{\frac{1}{n-1}}=\Sigma_{k,\beta_k}(B)\left(\frac{|\partial B|+\frac{\rho_n}{\eta_n}k^{\frac{n-2}{n(n-1)}}|B|}{k}\right)^{\frac{1}{n-1}}.
$$
Then, we can choose $\eps=\eps_k<\eps_0(\beta_k)$ sufficiently small so that
\begin{multline*}
\sigma_k(\Omega_{\beta_k,\eps_k})\left(\frac{|\partial\Omega_ {\beta_k,\eps_k}|}{k}\right)^{\frac{1}{n-1}}\geq \frac 12 \Sigma_{k,\beta_k}(B)\left(\frac{|\partial B|+\frac{\rho_n}{\eta_n}k^{\frac{n-2}{n(n-1)}}|B|}{k}\right)^{\frac{1}{n-1}}\\
\geq \frac 12\alpha_k\left(\frac{|\partial B|+\frac{\rho_n}{\eta_n}k^{\frac{n-2}{n(n-1)}}|B|}{k}\right)^{\frac{1}{n-1}}=\frac 12\eta_n\left(|\partial B|+\frac{\rho_n}{\eta_n}k^{\frac{n-2}{n(n-1)}}|B|\right)^{\frac{1}{n-1}}\geq C_n' k^{\frac{n-2}{n(n-1)^2}}.
\end{multline*}
for some $C_n'>0$ depending only on $n$. For $k\geq k_0(n)$ we set $\Omega_k:=\Omega_{\beta_k,\eps_k}$. For $1\leq k<k_0(n)$ we set $\Omega_k=B$. Since $\sigma_k(B)>0$ for all $k\geq 2$, after possibly changing the value of $C_n'$ inequality \eqref{lower_b} holds for every $k\geq 2$. The conclusion of Theorem \ref{main} follows since $n\geq 3$.

\smallskip

In order to conclude the proof, we show the validity of \eqref{final_claim}, that is:
$$
\Sigma_{k,\beta_k}(B)\geq\alpha_k.
$$ 

Assume by contradiction that  $\Sigma_{k,\beta_k}(B)<\alpha_k$. Now, the function $\alpha\mapsto\lambda_{k,\alpha}(B)$ is continuous and non-increasing in $\alpha$ (this follows from the min-max principle of Robin eigenvalues), hence, using Lemma \ref{lowerR} we deduce:
$$
S_n\beta_k\Sigma_{k,\beta_k}(B)=\lambda_{k,\Sigma_{k,\beta_k}(B)}(B)\geq\lambda_{k,\alpha_k}(B)\geq\rho_n k^{2/n}=S_n\beta_k\alpha_k
$$
hence
$$
\Sigma_{k,\beta_k}(B)\geq\alpha_k,
$$
a contradiction.
\qed

\section{Proof of Theorem \ref{main2}}\label{proof_main2}

We will prove Theorem \ref{main2} assuming $n\geq 3$, since for $n=2$ it follows directly from \eqref{isoperimetric2} and from the isoperimetric inequality in the plane: $|\partial\Omega|^2\geq4\pi|\Omega|$. Note that the quantity $|\Omega|^\frac 1n\sigma_k(\Omega)$ is scaling invariant, hence it is enough to prove the theorem assuming
$$
|\Omega|=1.
$$
Fix $k\geq 2$.

\smallskip

{\bf Case 1.} Assume that $|\partial\Omega|\geq k^\frac{n-2}{n(n-1)}$. We recall the isoperimetric bound \eqref{isoperimetric2} from \cite{colboisgirouard_steklov}:
\begin{equation}\label{isoperimetricbis}
|\partial\Omega|^\frac{1}{n-1}\sigma_k(\Omega)\leq\frac{C_n}{I(\Omega)^\frac{n-2}{n-1}}k^{\frac 2n}
\end{equation}
where $I(\Omega)=\frac{|\partial\Omega|}{|\Omega|^\frac{n-1}{n}}$. Since $|\Omega|=1$, from \eqref{isoperimetricbis} we get
$$
\sigma_k(\Omega)\leq C_n k^\frac{1}{n-1},
$$
which is the desired inequality.

\smallskip

{\bf Case 2.} Assume that $|\partial\Omega|\leq k^\frac{n-2}{n(n-1)}$. Thanks to the following Lemma, which we prove in Section \ref{23}, we will  build $k$ disjointly supported test functions with controlled Rayleigh quotient for $\sigma_k(\Omega)$. First, we define $Q_s(x)$ the cube of side $s$ centered at $x$:
$$
Q_s(x):=x+(-s/2,s/2)^n.
$$
Throughout the proof, $\delta_n,c_n$ will denote the following dimensional constants:
\begin{align*}
&\delta_n:=2^{-(2n+1)},\\
&c_n:=\delta_n^{n-1}=2^{-(2n+1)(n-1)}.
\end{align*}
\begin{lem}\label{lem_cap}
Let $0<a\leq\frac {1}{n-1}$ and let $\Omega$ be a smooth, bounded domain in $\mathbb R^n$, $n\geq 2$, with $|\Omega|=1$. If $k\geq 2$ and
\begin{equation}\label{small_per}
|\partial\Omega|\leq k^a,
\end{equation}
then there exist $k$ points $x_1,...,x_k\in\mathbb R^n$ such that the closed cubes 
$$
\overline Q_{2r}(x_i)
$$
with $r=\delta_n k^{-\frac{1}{n-1}}$ are pairwise disjoint and
\begin{equation}\label{area_decomp}
|\partial\Omega\cap Q_r(x_i)|\geq\frac{c_n}{k}
\end{equation}
for all $i=1,...,k$.
\end{lem}

We are now in position to conclude the proof of Theorem \ref{main2}. We apply Lemma \ref{lem_cap} with $a=\frac{n-2}{n(n-1)}$ (note that $a\leq\frac 1{n-1}$). We get points $x_1,...,x_k$ such that the closed cubes $\overline Q_{2r}(x_i)$, with $r=\delta_n k^{-\frac{1}{n-1}}$, are pairwise disjoint and satisfy $|\partial\Omega\cap Q_r(x_i)|\geq\frac{c_n}{k}$.

\smallskip For all $i=1,...,k$, let $\phi_i$ be a Lipschitz function supported in $Q_{2r}(x_i)$ such that $0\leq\phi_i\leq 1$,  $\phi_i\equiv 1$ in $Q_r(x_i)$ and $|\nabla\phi_i|\leq\frac 2r$ a.e. Then
$$
\int_\Omega|\nabla\phi_i|^2\leq\frac{4}{r^2}|Q_{2r}(x_i)|=2^{n+2}r^{n-2},
$$
while
$$
\int_{\partial\Omega}\phi_i^2\geq |\partial\Omega\cap Q_r(x_i)|\geq\frac{c_n}{k}.
$$
Hence
$$
\frac{\int_{\Omega}|\nabla\phi_i|^2}{\int_{\partial\Omega}\phi_i^2}\leq \frac{2^{n+2}}{c_n}kr^{n-2}=\frac{2^{n+2}}{\delta_n^{n-1}}k(\delta_nk^{-\frac{1}{n-1}})^{n-2}=\frac{2^{n+2}}{\delta_n}k^{\frac{1}{n-1}}.
$$
Let $W$ be the $k$-dimensional subspace of $H^1(\Omega)$ generated by $\phi_1,...,\phi_k$. Here $H^1(\Omega)$ is the usual Sobolev space of functions in $L^2(\Omega)$ with weak first derivatives in $L^2(\Omega)$. Since the supports of the $\phi_i$'s are disjoint, then the Rayleigh quotient of any $\phi\in W$ is bounded above by $\frac{2^{n+2}}{\delta_n}k^{\frac{1}{n-1}}$. The proof is concluded by recalling that
$$
\sigma_k(\Omega)=\min_{\substack{V\subset H^1(\Omega)\\{\rm dim}V=k}}\max_{0\ne\phi\in V}\frac{\int_{\Omega}|\nabla\phi|^2}{\int_{\partial\Omega}\phi^2}\leq \max_{0\ne\phi\in W}\frac{\int_{\Omega}|\nabla\phi|^2}{\int_{\partial\Omega}\phi^2}\leq\frac{2^{n+2}}{\delta_n}k^{\frac{1}{n-1}}=2^{3n+3}k^{\frac{1}{n-1}}.
$$
This concludes the proof of Theorem \ref{main2}.

\section{Proof of Lemmas \ref{lem_comparison}, \ref{countingR}, \ref{lowerR} and \ref{lem_cap}}\label{23}
\subsection{Proof of Lemma \ref{lem_comparison}}
We recall here the variational characterization of the eigenvalues of \eqref{dynamical} and \eqref{robin}:

\begin{equation}\label{minmaxdyn}
\Sigma_{k,\beta}(\Omega)=\min_{\substack{V\subset H^1(\Omega)\\{\rm dim}V=k}}\max_{0\ne u\in V}\mathcal R^{Dyn}_\beta(u),
\end{equation}
where
$$
\mathcal R^{Dyn}_\beta(u):=\frac{\int_\Omega|\nabla u|^2}{S_n\beta\int_\Omega u^2+\int_{\partial\Omega}u^2}
$$
and
\begin{equation}\label{minmaxrob}
\lambda_{k,\alpha}(\Omega)=\min_{\substack{V\subset H^1(\Omega)\\{\rm dim}V=k}}\max_{0\ne u\in V}\mathcal R^{Rob}_\alpha(u),
\end{equation}
where
$$
\mathcal R^{Rob}_\alpha(u):=\frac{\int_\Omega|\nabla u|^2-\alpha\int_{\partial\Omega}u^2}{\int_\Omega u^2}.
$$
Now, for any $u\in\ H^1(\Omega)$ we see that
$$
\mathcal R^{Rob}_\alpha(u)-\alpha S_n\beta=\frac{\int_{\partial\Omega}u^2+S_n\beta\int_\Omega u^2}{\int_\Omega u^2}\left(\mathcal R^{Dyn}_\beta(u)-\alpha\right).
$$
Hence
\begin{equation}\label{iff}
\mathcal R^{Rob}_\alpha(u)\leq\alpha S_n\beta\iff\mathcal R^{Dyn}_\beta(u)\leq\alpha.
\end{equation}
Now, set $\alpha=\Sigma_{k,\beta}(\Omega)$, and let $V_k\subset H^1(\Omega)$ be the subspace of $H^1(\Omega)$ spanned by the first $k$ eigenfunctions of the dynamical problem \eqref{dynamical}. From \eqref{iff} we deduce that
$$
\mathcal R^{Rob}_{\Sigma_{k,\beta}(\Omega)}(u)\leq \Sigma_{k,\beta}(\Omega) S_n\beta
$$
for all $u\in V_k\setminus\{0\}$, and from the min-max \eqref{minmaxrob}
$$
\lambda_{k,\Sigma_{k,\beta}(\Omega)}(\Omega)\leq S_n\beta\Sigma_{k,\beta}(\Omega).
$$

To obtain the reverse inequality, set $\alpha=\Sigma_{k,\beta}(\Omega)$ and let $W_k\subset H^1(\Omega)$ be the subspace of $H^1(\Omega)$ spanned by the first $k$ eigenfunctions of the Robin problem \eqref{robin}. From the min-max principle for the dynamical eigenvalues \eqref{minmaxdyn} we see that
$$
\Sigma_{k,\beta}(\Omega)\leq\max_{0\ne u\in W_k}\mathcal R^{Dyn}_\beta(u)
$$
and let $u^*\in W_k$ be such that
$$
\mathcal R^{Dyn}_\beta(u^*)=\max_{0\ne u\in W_k}\mathcal R^{Dyn}_\beta(u).
$$
From \eqref{iff} we get
$$
\mathcal R^{Rob}_{\Sigma_{k,\beta}(\Omega)}(u^*)\geq S_n\beta\Sigma_{k,\beta}(\Omega)
$$
but from the definition of $W_k$
$$
\mathcal R^{Rob}_{\Sigma_{k,\beta}(\Omega)}(u)\leq\lambda_{k,\Sigma_{k,\beta}(\Omega)}(\Omega)
$$
for all $u\in W_k$. We conclude that
$$
\lambda_{k,\Sigma_{k,\beta}(\Omega)}(\Omega)\geq S_n\beta\Sigma_{k,\beta}(\Omega).
$$
This concludes the proof. \qed

\subsection{Proof of Lemma \ref{countingR}}
As usual, we look for solutions to \eqref{robin} on $B$ of the form $u(r,\omega)=v(r)Y_\ell(\omega)$, where $(r,\omega)\in [0,1)\times\mathbb S^{n-1}$ and $Y_\ell$ is an eigenfunction of the Laplacian on $\mathbb S^{n-1}$ with eigenvalue $\ell(\ell+n-2)$, $\ell\in\mathbb N$. Hence the radial part $v$ solves
\begin{equation}\label{robinL}
\begin{cases}
-v''(r)-\frac{n-1}{r}v'(r)+\frac{\ell(\ell+n-2)}{r^2}v(r)=\lambda v(r)\,, & {\rm in\ }(0,1)\,,\\
v'(1)=\alpha v(1)\,,\\
\lim_{r\to 0}r^{n-1}v'(r)=0.
\end{cases}
\end{equation}
When $\ell=0$ the condition at $r=0$ is equivalent to $v'(0)=0$, while for $\ell\geq 1$ it is equivalent to $v(0)=0$. Now, for each $\ell\in\mathbb N$, problem \eqref{robinL} admits a sequence of simple eigenvalues:
$$
-\infty<\lambda_{1,\alpha}^{\ell}<\lambda_{2,\alpha}^{\ell}<\cdots<\lambda_{k,\alpha}^{\ell}<\cdots\nearrow+\infty
$$
and the spectrum of \eqref{robin} on $B$ is given by $\left\{\lambda_{k,\alpha}^{\ell}\right\}_{\ell\geq 0,k\geq 1}$, where each eigenvalue $\lambda_{k,\alpha}^\ell$ is counted with its multiplicity given by $d_{\ell,n}$, where  $d_{\ell,n}$ is the multiplicity of the Laplacian eigenvalue $\ell(\ell+n-2)$ on $\mathbb S^{n-1}$.

Let 
$$
H^1_\ell(0,1):=\{v\in H^1_{loc}(0,1):\int_0^1\left(v'^2+v^2+\frac{\ell(\ell+n-2)}{r^2}v^2\right)r^{n-1}dr<\infty\}.
$$
The eigenvalues of \eqref{robinL} are characterized by
\begin{equation}\label{minmaxR}
\lambda_{k,\alpha}^\ell=\min_{\substack{V\subset H^1_\ell(0,1)\\{\rm dim}V=k}}\max_{0\ne v\in V}\frac{\int_0^1\left(v'^2+\frac{\ell(\ell+n-2)}{r^2}v^2\right)r^{n-1}dr-\alpha v^2(1)}{\int_0^1v^2r^{n-1}dr}.
\end{equation}
Consider now the following family of one-dimensional Dirichlet problems indexed by $\ell\in\mathbb N$:

\begin{equation}\label{dirichletL}
\begin{cases}
-v''(r)-\frac{n-1}{r}v'(r)+\frac{\ell(\ell+n-2)}{r^2}v(r)=\lambda^D v(r)\,, & {\rm in\ }(0,1)\,,\\
v(1)=0\,,\\
\lim_{r\to 0}r^{n-1}v'(r)=0.
\end{cases}
\end{equation}
Now, each problem \eqref{dirichletL} admits a sequence of simple eigenvalues:
$$
0<\lambda_1^{\ell,D}<\lambda_{2}^{\ell,D}<\cdots<\lambda_{k}^{\ell,D}<\cdots\nearrow+\infty.
$$
Then, the set $\left\{\lambda_{k}^{\ell,D}\right\}_{\ell\geq 0,k\geq 1}$ (each counted with its multiplicity $d_{\ell,n}$) is the set of the eigenvalues of the Dirichlet Laplacian on $B$.

Let 
$$
H^1_{0,\ell}(0,1):=\{v\in H^1_\ell(0,1):v(1)=0\}.
$$
Hence $H^1_{0,\ell}(0,1)\subset H^1_\ell(0,1)$ has codimension $1$.
The eigenvalues of \eqref{dirichletL} are characterized by
\begin{equation}\label{minmaxD}
\lambda_k^{\ell,D}=\min_{\substack{V\subset H^1_{0,\ell}(0,1)\\{\rm dim}V=k}}\max_{0\ne v\in V}\frac{\int_0^1\left(v'^2+\frac{\ell(\ell+n-2)}{r^2}v^2\right)r^{n-1}dr}{\int_0^1v^2r^{n-1}dr}.
\end{equation}

From the min-max principle for the eigenvalues  of \eqref{robinL} and \eqref{dirichletL} and  from the fact that the codimension of $H^1_{0,\ell}(0,1)$ is $1$, we have that, for all $\alpha\in\mathbb R$, $\ell\in\mathbb N$ and $k\geq 1$:
\begin{equation}
\lambda_{k,\alpha}^\ell\leq\lambda_k^{\ell,D}\leq\lambda_{k+1,\alpha}^\ell,
\end{equation}
which implies
$$
N_\alpha^{\ell,R}(E)\leq N^{\ell,D}(E)+\#\{\lambda_{1,\alpha}^\ell< E\}
$$
where
$$
N_\alpha^{\ell,R}(E):=\#\{\lambda_{j,\alpha}^\ell<E\}\,,\quad\quad\quad N^{\ell,D}(E):=\#\{\lambda_{j}^{\ell,D}<E\}.
$$
Hence 
\begin{equation}\label{sum2}
N^R_\alpha(E)\leq N^D(E)+\sum_{\ell:\lambda_{1,\alpha}^\ell<E}d_{\ell,n},
\end{equation}
where $N^D(E)$ is the counting function for the Dirichlet eigenvalues on $B$.

 Now, let $u$ be an eigenfunction of \eqref{robin} with eigenvalue $\lambda_{1,\alpha}^\ell$. Integrating by parts we get 
 $$
\int_B\langle ux,\nabla u\rangle=\int_{\partial B}u^2-\int_B{\rm div}(ux)u=\int_{\partial B}u^2-\int_B\langle ux,\nabla u\rangle-n\int_B u^2,
 $$
 which implies, when $\alpha>0$
\begin{equation}\label{A}
\int_{\partial B}u^2=n\int_Bu^2+2\int_B\langle ux,\nabla u\rangle\leq\frac{1}{2\alpha}\int_B|\nabla u|^2+(n+2\alpha)\int_Bu^2.
\end{equation}
 Moreover, since $u=v(r)Y_\ell(\omega)$, we have, for all $\alpha\geq 0$
 \begin{multline}\label{B}
\int_B|\nabla u|^2=\int_0^1\int_{\mathbb S^{n-1}}\left(v'(r)^2Y_\ell^2(\omega)+\frac{v^2(r)|\nabla Y_\ell|^2}{r^2}\right)r^{n-1}dv_{\mathbb S^{n-1}}dr\\
\geq\int_0^1\int_{\mathbb S^{n-1}}v^2(r)|\nabla Y_\ell|^2r^{n-3}dv_{\mathbb S^{n-1}}dr\\
=\ell(\ell+n-2)\int_0^1\int_{\mathbb S^{n-1}}v^2(r)Y^2_\ell(\omega)r^{n-3}dv_{\mathbb S^{n-1}}dr\\\geq\ell(\ell+n-2)\int_0^1\int_{\mathbb S^{n-1}}v^2(r)Y^2_\ell(\omega)r^{n-1}dv_{\mathbb S^{n-1}}dr
=\ell(\ell+n-2)\int_Bu^2
 \end{multline}
 From \eqref{A}, \eqref{B} and the min-max principle for $\lambda_{1,\alpha}^\ell$ we conclude that
 $$
\lambda_{1,\alpha}^\ell\geq\frac{1}{2}\ell(\ell+n-2)-(n\alpha+2\alpha^2).
 $$
 This conclusion follows just from \eqref{B} if $\alpha=0$.
 Hence 
 $$
\lambda_{1,\alpha}^\ell<E \Rightarrow \ell\leq C_{n,1}(1+\alpha+\sqrt{E})
 $$
 where $C_{n,1}$ depends only on $n$. Now, since
 $$
d_{\ell,n}=\frac{(2\ell+n-2)(\ell+n-3)!}{\ell!(n-2)!}
$$
 we conclude that
 \begin{equation}\label{sum_mult}
\sum_{\ell:\lambda_{1,\alpha}^\ell<E}d_{\ell,n}\leq C_{n,2} (1+\alpha+\sqrt{E})^{n-1},
 \end{equation}
 where $C_{n,2}$ depends only on $n$.
  From the well-known Berezin-Li-Yau upper bound on $N^D(E)$ \cite{Berezin,LiYau}:
\begin{equation}\label{BLY}
N^D(E)\leq C_{n,3}(1+E)^{n/2},
  \end{equation}
where $C_{n,3}$ depends only on $n$. The proof is concluded by using \eqref{sum_mult} and \eqref{BLY} in \eqref{sum2}. \qed

\subsection{Proof of Lemma \ref{lowerR}}
Let $\rho_n$, $\eta_n$ and $k_0(n)$ and $\alpha_k$ be as in the statement of Lemma \ref{lowerR}:
\begin{eqnarray*}
&&\rho_n:=(8(1+C_n))^{-2/n}\,,\\
&&\eta_n:=\frac{1}{2}(8(1+C_n))^{-1/(n-1)}\,,\\
&&k_0(n):=\left[2^{4n(n-1)}(1+C_n)\right]+1\,,\\
&&\alpha_k:=\eta_n k^{1/(n-1)}.
\end{eqnarray*}
To simplify the notation, throughout the proof we write $\rho,\eta$ in place of $\rho_n$ and $\eta_n$. Define, for any $k\in\mathbb N$, $k\geq 1$:
$$
E_k=\rho k^{2/n}.
$$
From Lemma \ref{countingR}:
$$
N_{\alpha_k}^R(E_k)\leq C_n\left(1+\rho^{n/2}k+(\eta^2k^{2/(n-1)}+\rho k^{2/n})^{(n-1)/2}\right),
$$
from which
\begin{equation}\label{partialA}
\frac{N_{\alpha_k}^R(E_k)}{k}\leq \frac{C_n}{k}+C_n\rho^{n/2}+C_n (\eta^2+\rho k^{-2/(n(n-1))})^{(n-1)/2}
\end{equation}
We note that, from the definition of $\rho$ we have $C_n\rho^{n/2}\leq\frac 18$.
Now, for $k\geq k_0(n)$ we have:
$$
\frac{C_n}{k}\leq 2^{-4n(n-1)}\leq\frac{1}{16}.
$$
and
$$
\frac{\rho k^{-2/(n(n-1))}}{\eta^2}=4(8(C_n+1)/k)^{2/(n(n-1))}<1.
$$
Hence
$$
\eta^2+\rho k^{-2/(n(n-1))}\leq 2\eta^2
$$
which implies
\begin{equation}\label{partialB}
C_n(\eta^2+\rho k^{-2/(n(n-1))})^{(n-1)/2}\leq C_n 2^{(n-1)/2}\eta^{n-1}\leq\frac 1{16}
\end{equation}
since $n\geq 3$. Using \eqref{partialB} in \eqref{partialA} we get
$$
\frac{N_{\alpha_k}^R(E_k)}{k}\leq\frac 14<1
$$
from which
$$
\lambda_{k,\alpha_k}(B)\geq \rho k^{2/n}.
$$
This concludes the proof.

\qed

\subsection{Proof of Lemma \ref{lem_cap}}
Define 
$$
w(x):=|\Omega\cap Q_r(x)|\,,\ \ \ E:=\{x\in\mathbb R^n:w(x)>\frac{r^n}{2}\}\,,\ \ \ D:=\Omega\setminus E.
$$
The function $w$ is continuous hence $E$ is open and bounded. The proof is divided in three steps.

\smallskip

{\bf Step 1.} Here we prove that $E$ keeps a quite large part of the volume of $\Omega$. In particular, we prove that
$$
|E|\geq\frac 34.
$$
If $D=\emptyset$, there is nothing to prove. Otherwise, we choose a maximal $r/2$-separated family of points $p_1,...,p_N\in D$ with respect to the norm $|\cdot|_\infty$. We recall that, for points $x=(x_1,...,x_n)$, $y=(y_1,...,y_n)$ in $\mathbb R^n$,  $|x-y|_\infty=\max_{i=1,...,n}|x_i-y_i|$. Here by maximality we mean that we cannot add further points of $D$ to the family $p_1,...,p_N$ and keep $|p_i-p_j|_\infty\geq r/2$ (when $i\ne j$). This family exists since $D$ is bounded. Also, by maximality
\begin{equation}\label{inclusionD}
D\subset\bigcup_{i=1}^NQ_r(p_i).
\end{equation}
Moreover, no point $y$ belongs to more than $4^n$ cubes $Q_r(p_i)$, that is
\begin{equation}\label{boundN}
\#\{i:y\in Q_r(p_i)\}\leq 4^n.
\end{equation}
In fact, if $y\in Q_r(p_i)$ for $i=1,...,m$ for some $m$ (possibly after re-labeling the cube indices), then the pairwise disjoint cubes $Q_{r/2}(p_i)$, $i=1,..,m$, would be contained in $Q_{2r}(y)$. Hence $m(r/2)^n\leq (2r)^n$, so that $m\leq 4^n$.

\smallskip

We recall the following sharp {\it relative isoperimetric inequality} in cubes, originally proved by Hadwiger for polyhedral sets \cite{hadwiger}, see also \cite[Formula (3.20)]{ABBF}: if $A$ is a smooth open set in $\mathbb R^n$ and $Q_r$ is a cube of side $r$ 
\begin{equation}\label{relativeiso}
\min\{|A\cap Q_r|,r^n-|A\cap Q_r|\}\leq\frac{r}{2}|\partial A\cap Q_r|.
\end{equation}
We  apply this inequality to  $A=\Omega$ and $Q_r=Q_r(p_i)$. Since $p_i\in D$, by definition of $D$, we have
$$
|\Omega\cap Q_r(p_i)|\leq\frac{r^n}{2},
$$
and hence, from \eqref{relativeiso}:
$$
|\Omega\cap Q_r(p_i)|\leq\frac{r}{2}|\partial\Omega\cap Q_r(p_i)|.
$$
From \eqref{inclusionD} and \eqref{boundN} we find that
$$
|D|\leq 2^{2n-1}r|\partial\Omega|\leq\frac{1}{4},
$$
where the last inequality follows since $0<a\leq\frac{1}{n-1}$ and $\delta_n=2^{-(2n+1)}$. Then, since $|\Omega|=1$, we have
$$
|E|\geq\frac{3}{4}, 
$$
which is the claim.

\smallskip

{\bf Step 2.} Here we prove that $\partial E$ contains at least $k$ well-separated points. $E$ is non-empty, open and bounded, hence $\partial E$ is compact and we can choose a maximal, $3r$-separated family of points $q_1,...,q_M\in\partial E$ with respect to the norm $|\cdot|_\infty$. By maximality
\begin{equation}\label{pdEinclusion}
\partial E\subset\bigcup_{i=1}^M Q_{6r}(q_i).
\end{equation}
Now, for a point $x=(x_1,...,x_n)\in\mathbb R^n$ we define
$$
\pi_j(x):=(x_1,...,x_{j-1},x_{j+1},...,x_n)
$$
its projection on the orthogonal hyperplane to the $j$-th coordinate direction. We have the elementary inclusion
\begin{equation}\label{projection}
\pi_j(E)\subset\pi_j(\partial E)
\end{equation}
for all $j=1,...,n$. From \eqref{pdEinclusion} and \eqref{projection} we get
$$
|\pi_j(E)|\leq M (6r)^{n-1},
$$
and by the fact that $|A|^{n-1}\leq\Pi_{j=1}^n|\pi_j(A)|$ for any bounded, open set $A$ in $\mathbb R^n$ (see e.g., \cite{loomis_w}), we conclude
$$
|E|^{n-1}\leq\Pi_{j=1}^n|\pi_j(E)|\leq M^n(6r)^{n(n-1)},
$$
which, together with $|E|\geq\frac 34$ proved in Step 1, $6\delta_n=6\cdot 2^{-(2n+1)}\leq\frac 12$ and $r=\delta_n k^{-\frac 1{n-1}}$, implies
$$
M\geq k.
$$

{\bf Step 3.} We conclude the proof. Choose $k$ of the points $q_1,...,q_M$ found in Step 2, and relabel them as $x_1,...,x_k$. Since $q_i$ are $3r$-separated, then the closed cubes $\overline Q_{2r}(x_i)$ are pairwise disjoint.

Recall that $E=\{x:w(x)>r^n/2\}$, and $w(x)=|\Omega\cap Q_r(x)|$ is continuous. Since $x_i\in\partial E$:
$$
\frac{r^n}{2}=|\Omega\cap Q_r(x_i)|\leq\frac{r}{2}|\partial\Omega\cap Q_r(x_i)|,
$$
where we have used again the relative isoperimetric inequality \cite[Formula (3.20)]{ABBF}. This is equivalent to
$$
|\partial\Omega\cap Q_r(x_i)|\geq r^{n-1}=\frac{\delta_n^{n-1}}{k}=\frac{c_n}{k}.
$$
The proof is concluded. \qed

\begin{rem}\label{rem:riemannian}
The only genuinely Euclidean ingredient of the above proof is Step 2. Steps 1 and 3
are local and survive on a fixed Riemannian manifold, possibly with geodesic balls in place of
cubes and different constants. On the other hand, the global argument in Step 2 based on the coordinate projections, the
inclusion $\pi_j(E)\subset\pi_j(\partial E)$ and the Loomis--Whitney inequality have
no analogue on a closed manifold. In fact, let $\Omega_{\varepsilon}$ be the
 domains of \cite{GiLa} obtained, for example, by periodically perforating the round sphere $\mathbb S^n$  (with constant weight $\beta$) and, for $k$ sufficiently large, let $r=\delta_n k^{-\frac{1}{n-1}}$. Then, for $\eps$ sufficiently small (depending on $k$), we have $|B_r(p)\cap\Omega_\eps|>|B_r(p)|/2$ everywhere, so that the corresponding set $E$ is the whole manifold and $\partial E=\emptyset$,
while $|\partial\Omega_{\varepsilon}\cap B_r(p)|\leq c r^{n}\ll r^{n-1}= c_n/k$,
uniformly in $p\in \mathbb S^n$ (see also \cite[Proof of Lemma 3.2]{GiLa}). Here $B_r(p)$ denotes the geodesic ball of radius $r$ centered at $p$. Thus the analogue of conclusion \eqref{area_decomp} of Lemma \ref{lem_cap} fails in this case, as it should, otherwise the test functions argument in the proof of Theorem \ref{main2} would imply the wrong $k^\frac{1}{n-1}$ term in the volume-normalised upper bound for the family $\Omega_\eps$.
\end{rem}

\section*{Acknowledgments}
The author is grateful to Bruno Colbois for first introducing him to the study of upper bounds for Steklov eigenvalues in 2016, and to Joachim Stubbe for proposing to pursue this study further together. The author is also grateful to Alexandre Girouard and Jean Lagacé for helpful comments and stimulating exchanges concerning homogenisation methods and the question of the optimal eigenvalue growth rate.

\section*{AI usage disclosure}
The author acknowledges the use of ChatGPT 5.6 Sol for mathematical discussions and
editorial assistance. In particular, the AI tool contributed to the development of the proofs of  Lemmas \ref{lem_comparison}, \ref{countingR}, \ref{lowerR} and \ref{lem_cap}.  All AI-assisted arguments and computations were independently checked by the author, who takes full responsibility for the contents of the paper.

\bibliography{bibliography.bib}
\bibliographystyle{abbrv}

\end{document}